\documentclass{amsart}[12pt,a4paper]

\usepackage[applemac]{inputenc}
\usepackage[colorlinks=true,pagebackref]{hyperref}
\usepackage[english]{babel} 
\usepackage{color}  
\usepackage{latexsym}
\usepackage[T1]{fontenc} 
\usepackage{graphicx}
\usepackage{fancyhdr}
\usepackage{enumerate}

\usepackage{amsfonts}  
\usepackage{amsmath}
\usepackage{amssymb} 
\usepackage{amsthm}
\input xy
\xyoption{all}

\DeclareRobustCommand{\SkipTocEntry}[5]{}

\newtheorem{theorem}{Theorem}[section]
\newtheorem{proposition}[theorem]{Proposition}
\newtheorem{lemma}[theorem]{Lemma}
\newtheorem{definition}[theorem]{Definition}
\newtheorem{corollary}[theorem]{Corollary}
\newtheorem{def.thm}[theorem]{Definition-Theorem}
\newtheorem{question}[theorem]{Question}

\theoremstyle{definition}
\newtheorem{remark}[theorem]{Remark}

\newcommand{\defn}{:=}

\newcommand{\natl}{\mathbb{N}}

\newcommand{\real}{\mathbb{R}}
\newcommand{\complex}{\mathbb{C}}
\newcommand{\ratl}{\mathbb{Q}}
\DeclareMathOperator{\NE}{\overline{NE}_1}

\newcommand{\strutt}{\mathcal{O}}
\newcommand{\nklt}{\mathrm{Nklt}}
\newcommand{\nlc}{\mathrm{Nlc}}

\newcommand{\cdel}{\Delta}

\newcommand*{\QEDB}{\hfill\ensuremath{\square}}%

\begin{document}

\title{Hyperbolicity for log canonical pairs and the cone theorem}

\author{Roberto Svaldi}
\date{\today}	
\address{Roberto Svaldi, EPFL SB MATH MATH-GE,
	MA B1 497 (B\^{a}timent MA), Station 8,
	CH-1015 Lausanne, Switzerland.}
\email{roberto.svaldi@epfl.ch}


\begin{abstract}
\textsl{Given a log canonical pair $(X, \Delta)$,
we show that $K_X+\Delta$ is nef assuming there is no non-constant map from the affine line 
with values in the open strata of the stratification induced by the non-klt locus of $(X, \Delta)$. 
This implies a generalization of the Cone Theorem where each $K_X+\Delta$-negative extremal ray is spanned by a rational curve that is the closure of a copy of the affine line contained in one of the open strata of $\nklt(X, \Delta)$.
Moreover, we give a criterion of Nakai type to  determine when under the above condition $K_X+\Delta$ is ample and we prove some partial results in the case of arbitrary singularities.}
\end{abstract}

\bibliographystyle{amsalpha}

\maketitle

\tableofcontents

\section{Introduction}

Understanding the existence and distribution of curves on a given algebraic variety is a classical problem in algebraic geometry. 
For example, its significance in understanding the birational structure of algebraic varieties geometry -- in particular, for the case of rational curves -- has been evident since the early days of the subject, when the Italian School started the classification of algebraic surfaces.

In the past 30 years, with the emergence and development of the so-called Minimal Model Program (in short, MMP) such aspect has been investigated and understood in birational geometry in far greater generality, thanks to the work of many different people. 
The main realization has been that the existence of rational curves on a mildly singular normal variety is strictly related to the positivity properties of the cotangent bundle.

On the other hand, rational curves on varieties have been object of study long before the MMP was even imagined. 
Many authors turned their attention to the study of the existence/absence of rational curves and their distribution on a given variety, providing some interesting discoveries and conjectures. 
There are a number of famous open questions due to several authors that similarly predict a strong link between the positivity of the curvature of the cotangent bundle of a variety $X$ and the absence or bounded distribution  of non-trivial holomorphic maps $f \colon \mathbb{C} \to X$. 
The interested reader can consult \cite{dem2012} for a survey of classical and more recent questions and
results in this direction.

The main result of this paper is inscribed in this line of thought: we show that there is a clear connection between positivity properties of pairs given by algebraic varieties together with an effective divisor and the hyperbolicity of a stratification that is naturally induced by the singularities of the divisor and the ambient variety.

\begin{theorem}
\label{main.thm}
Let $X$ be a smooth projective variety and $D = \sum_{j \in J} D_j$ be a reduced simple normal crossing divisor on $X$. 
Assume that 
  \begin{itemize}
   \item there is no non-constant morphism $f : \mathbb{A}^1 \to X \setminus D$
   \item for any intersection of components of $D, \; D_I =\cap_{i \in I} D_i, \;  I \subset J$ there is no non-constant morphism $f : \mathbb{A}^1 \to (D_I \setminus \cup_{j \in (J \setminus I)} D_j)$.
  \end{itemize}
Then $K_X + D$ is nef.
\newline
More generally, let $(X, \Delta)$ be a log canonical pair. 
Assume that there is no non-constant morphism 
\[
f \colon \mathbb{A}^1 \to X \setminus \{x \in X \; | \; (X, \Delta) \textrm{ is not Kawamata log terminal at } x\}
\] 
and the same holds for all the open strata of the non-klt locus.
\newline
Then $K_X + \Delta$ is nef.
\end{theorem}

Following Lu and Zhang, {\cite{lu-zh}}, we say that a pair $(X, \Delta)$ is Mori hyperbolic if it satisfies the assumptions in the above theorem on the non-existence of copies of $\mathbb{A}^1$ in the (open) stratification induced by the non-klt locus of $(X, \Delta)$.
We generalize this definition to any normal singularity in Definition \ref{m.hyp.def}. 
When $(X, \Delta)$ is a simple normal crossing pair, then the non-klt locus of $\Delta$, denoted $ \nklt(X, \Delta)$, is the union of the components of coefficients $\geq$ 1 in $\Delta$. 
If the pair is not simple normal crossing, the non-klt locus is the image of the components of coefficient $\geq 1$ of the pullback of $K_X + \Delta$ to a log resolution, cf. Section \ref{nonlc.locus.sec}.

Lu and Zhang proved a version of Theorem \ref{main.thm} for divisorial log terminal pairs assuming some factoriality conditions on the components, {\cite[Thm. 3.1]{lu-zh}}. 
Similar results, in the context of algebraic stacks -- and hence coarse moduli with quotient singularities -- were obtained by McQuillan and Pacienza in {\cite{mcq-pac}}. 
Theorem \ref{main.thm} reproves these results and moreover shows it can be extended to the category of log canonical pairs: this is a much larger class of pairs that do not necessarily have rational singularities. 
Hence, a priori it is not clear why the stratification described in the statement of the theorem should contain rational curves at all.
Theorem \ref{main.thm} has also been extended to the category of rank 2 foliations on algebraic threefolds in \cite{mio.term}

The starting point for the MMP in the 1980s, was the discovery, due to Mori -- later improved by Koll\'ar, Reid, Shokurov, Kawamata, Ambro -- that the portion of the effective cone of curves on a normal mildly singular variety $X$ generated by classes of negative intersection with the canonical divisor $K_X$  is actually spanned by countably many classes of rational
curves. 
This is a now classical result that goes under the name of Cone Theorem, cf. {\cite[Thm. 1.24]{koll-mor-book}}. It has been generalized to divisors of the form $K_X + \Delta$, when the pair $(X, \Delta)$ has suitably
nice singularities. 
It immediately implies that the absence of rational curves on a variety $X$ guarantees the nefness
of $K_X + \Delta$. 
Nonetheless, that is an extremely strong assumption.
In order to obtain statements that apply to a wider class of cases,  one is lead to wonder what kind of hyperbolicity-like assumptions a  pair $(X, \Delta)$ could satisfy for $K_X + \Delta$ to be nef. 
Moreover, in case such assumptions are not satisfied,  one could then try to investigate how rational curves are distributed with respect to $\Delta$.

For example, let us consider a smooth quasi-projective variety $U$ and a compactifying simple normal crossing pair $(X, \Delta)$, $U= X \setminus \Delta$. 
In such context, these questions make even more sense in view of Iitaka's principle, see \cite[pg. 112]{matsuki.book}.
In this context, Iitaka's principle is just predicting a correspondence between theorems about non-singular varieties and regular differential forms and theorems about quasi-projective varieties and their regular differential forms which extend to the boundary of a compactification with at worst poles of order 1.

Using Theorem \ref{main.thm}, we are able to establish a version of the Cone Theorem describing the distribution of rational curves spanning $(K_X + \Delta)$-negative extremal rays with respect to the boundary $\Delta$. 
The Cone Theorem for lc pairs, \cite{amb03, fuj.nonvan.12} is a natural extension of the classical version of the cone theorem, \cite[Thm. 3.7]{koll-mor-book}.
Part of the difficulty lies in controlling what happens along the non-klt locus, as that is the most delicate locus where some of the classical results of the MMP may fail.
Using the techniques of Theorem \ref{main.thm}, we can improve the classical version of the Cone Theorem by providing precise information about the position of the rational curves that span the $(K_X+\Delta)$-negative extremal rays, with respect to the non-klt locus of $(X, \Delta)$.
More precisely, we can show that copies of either the affine or the projective naturally appear in the open strata of the stratification induced by lc centers on the non-klt locus.

We shall use $\NE(X)$ to denote the closure of the cone spanned by effective curves inside the group of curves with real coefficients modulo numerical equivalence.

\begin{theorem}[= Theorem~\ref{generalized.cone.thm2}]
  \label{generalized.cone.thm}
Let $(X, \Delta)$ be a log canonical pair.
\newline  
There exists countably many $K_X + \Delta$-negative rational curves $C_i$ such that
\[ 
\NE (X) = \NE (X)_{K_X + \Delta \geq 0} + \sum_{i\in I} \mathbb{R}_{>0}[C_i]. 
\]
Moreover, one of the two following conditions hold:
\begin{itemize}
\item $C_i \cap (X \setminus  \nklt(X, \Delta))$ contains the image of a non-constant morphism $f \colon \mathbb{A}^1 \to X$;
    
\item there exists an open stratum $\overline{W}$ of $\nklt(X, \Delta)$ such that $C_i \cap W$ contains the image of a non-constant morphism $f \colon \mathbb{A}^1 \to W$.
\end{itemize}
\end{theorem}

When $(X, \Delta)$ is a simple normal crossing pair, then the appearance of morphisms $f \colon \mathbb{A}^1 \to X \setminus D$ should be thought as the realization the Iitaka principle for the Cone Theorem. 

Finally, for a Mori hyperbolic pair $(X, \Delta)$, we prove that the classical Nakai-Moishezon-Kleiman criterion, {\cite[Thm. 1.2.23]{laz1.book}}, can be restated in a much simpler form: namely, it is enough to test ampleness only along the (finitely many) lc centers of $\Delta$ rather than having to check positivity of the self-intersection numbers along all subvarieties of $X$.

\begin{theorem}[= Corollary~\ref{nakai.crit.m.hyp.cor}]
\label{ample.thm}
Let $(X, \Delta)$ be a dlt pair. Assume that the pair is Mori hyperbolic.
\newline  
  Then the following are equivalent:
  \begin{enumerate}
    \item[(i)] $K_X + \Delta$ is ample;
    
    \item[(ii)] $(K_X + \Delta)^{\dim X} > 0$ and $(K_X + \Delta)^{\dim W} \cdot W > 0$ for every log canonical center $W \subset X$ of $(X, \Delta)$.
  \end{enumerate}
\end{theorem}

\subsection*{Sketch of the proof}
We explain now the structure of the proof of Theorem \ref{main.thm}.

The notion of Mori hyperbolicity for a log pair $(X, \Delta)$ has an inherently inductive nature. 
Hence, it is fair to expect that some sort of inductive approach could possibly lead to the above theorem. 
Indeed, this is the philosophy that we adopt in the course of the proof. 
A fundamental step in this sense is represented by the following result which makes clear the connection between the positivity of a Mori hyperbolic pair and its positivity along the non-klt locus of $(X, \Delta)$. That is in fact a general guiding principle in the study of purely lc pairs.

\begin{theorem}
  \label{nef.induct.cor}
  Let $(X, \Delta)$ be a log pair. 
  Assume that $(X, \Delta)$ is Mori hyperbolic. 
  Then $K_X + \Delta$ is nef if it is nef when restricted to the non-klt locus of $(X, \Delta)$.
\end{theorem}

To be able to use this result, we are actually forced to deal with singularities worse than log canonical. 
In the simple normal crossing case, in fact, Theorem \ref{nef.induct.cor} immediately implies Theorem \ref{main.thm} simply by performing adjunction along the components of $\Delta$ of coefficient $1$ and by using Kawamata's estimates on the length of extremal rays, {\cite{kawam91}}. 

In the log canonical case, instead, the strata of the non-klt locus of $(X,\Delta)$ are not as well behaved as in the simple normal crossing. 
It is just not possible to perform adjunction along a divisorial component, as there may not be any. 
Because of this, one tries to construct a new log pair $(X', \Delta')$ with positive coefficients and nice singularities (of dlt type) together with a birational morphism $\pi \colon  X' \to X$ such that $K_{X'} + \Delta' = \pi^{\ast} (K_X + \Delta)$, cf. Theorem \ref{dlt.mod.thm}. 
The proof is then carried out by conducting a careful analysis of adjunction along lc centers of codimension greater than 1 with respect to the morphism $\pi$, by means of the canonical bundle formula. 
It is in the course of this last part of the proof that we have to consider also log pairs with singularities worse then log canonical. 
This is the truly new insight which is needed to generalize the whole result to the log canonical case and for which Theorem \ref{nef.induct.cor} has been developed.

The paper is structured as follows: in Section \ref{sing.mmp.sect} and \ref{subadj.sect}, we recall some preliminaries about singularities of the Minimal Model Program and adjunction theory for lc centers of codimension higher than one. 
In Section \ref{dlt.mod.sect.}, we prove a special version of the existence of dlt modifications that will be needed in the proof of the Theorem \ref{main.thm}. 
In Section \ref{mh.sect}, we define Mori hyperbolicity and describe some of its properties. 
Section \ref{proof.sect} is devoted to the proof of Theorems \ref{main.thm} and \ref{generalized.cone.thm},  while in Section \ref{ample.sect} we prove Theorem \ref{ample.thm}.

\addtocontents{toc}{\SkipTocEntry}

\subsection*{Acknowledgments} The author wants to thank his Ph.D.advisor, 
Professor James M\textsuperscript{c}Kernan, for offering constant support and
encouragement during the development of this project. 
He thanks John Lesieutre and Tiankai Liu for useful discussions and
exchange of ideas, as well as Professors Gabriele di Cerbo, Antonio Lerario, and Jacopo Stoppa for reading 
preliminary versions of this work and for their 
encouragement and Prof. Steven Lu for his interest in this work.
He also wishes to thank the anonymous referees for the beneficial suggestions related to the structure of the paper.

Part of this work was done while the author was supported by NSF DMS \#0701101 and \#1200656.
The author would also like to thank MIT where he was a graduate student when most of this 
work was completed. During the final revision of this work the author was visiting Professor Jacopo Stoppa at SISSA, Trieste. 
He would to thank Jacopo Stoppa and SISSA for providing such a congenial place to work. 
The visit was supported by funding from the European Research Council under the European Union's Seventh Framework 
Programme (FP7/2007-2013)/ERC Grant agreement no. 307119.

This work is part of the author's Ph.D. thesis which received the 2016 Federigo Enriques Prize, awarded by the Unione Matematica Italiana and Centro Studi Enriques. 
The author wishes to thank these institutions for such honor.

\addtocontents{toc}{\SkipTocEntry}
\subsection*{Notation and Conventions}

By the term variety, we will always mean an integral, separated, projective scheme over an algebraically closed field $k$. 
Unless otherwise stated, it will be understood that $k = \complex$.
\newline
Unless otherwise specified, we adopt the same notations and conventions as in \cite{koll-mor-book}.
\\

\noindent 
If $D=\sum d_iD_i$ is an $\mathbb{R}$-divisor on a normal variety $X$, where the $D_i$'s are the distinct prime components of $D$, then we define $D^{*c} \defn \sum_{d_i \; * \; c} d_i D_i, \; c \in \real$, where $*$ is any of $=, \geq, \leq, >, <$.
\newline
The \textit{support} of an $\mathbb{R}$-divisor $\Delta=\sum_{i\in I}d_iD_i$ is the union of the prime divisors appearing in the formal sum, $\mathrm{Supp}(D)=\bigcup_{\{i\in I  \; | \; d_i\neq0\}}D_i$.
\newline
A \textit{sub-log pair} $(X,\Delta)$ consists of a normal variety $X$ and a Weil$\mathbb{R}$-divisor $\Delta$ such that $K_X+\Delta$ is $\mathbb{R}$-Cartier. 
If $\Delta$ is effective then we say that the sub-log pair $(X, \Delta)$ is a {\it log pair}.
\newline
A log pair $(X, \Delta= \sum_{i \in I} D_i)$, where $\Delta$ is an effective reduced divisor, is simple normal crossing (in short, snc) if $X$ as well as every component of $\Delta$ are smooth and, moreover, all components $D_i$ of $D$ intersect as transversally as possible, i.e. for every $p \in X$ one can choose a neighborhood $U \ni p$ (in the Zariski topology) and local coordinates $x_j$ s.t. for every $i$ there is an index $c(i)$ for which $D_i \cap U = (x_{c(i)} = 0)$. 
If $(X,\Delta)$ is snc a \textit{stratum} of $(X,\Delta)$ is either $X$ or an irreducible component of the intersection $\cap_{\{i\in I \; | \; d_i=1\}}D_j$.
Given a (closed) stratum, $W$, the corresponding \textit{open stratum} is obtained from $W$ by removing all the strata contained in $W$.
\newline
Given a normal variety $X$, a $\mathbb{K}$-b-divisor is a (possibly infinite) sum of geometric valuations of $k(X)$ with coefficients in $\mathbb{K}$,
\[
 \mathbb{D}= \sum_{i \in I} b_i V_i, \; V_i \subset k(X) \; {\rm and} \; b_i \in \mathbb{K}, \; \forall i \in I, \; 
\]
such that for every normal variety $X'$ birational to $X$, only a finite number of the $V_i$ can be realized by divisors on $X'$. 
The trace of $\mathbb{D}$ on $X'$, $\mathbb{D}_{X'}$, is defined as 
\[
\mathbb{D}_{X'}= \sum_{\substack{i \in I \\ c_{X'}(V_i)= D_i \\ D_i \; {\rm is \; a\; divisor}}} b_i D_i.
\]

\section{Pairs and their singularities}\label{sing.mmp.sect}

\begin{definition}

A \textit{log resolution} for a sub-log pair $(X,\Delta)$ is a projective birational morphism $\pi \colon X' \to X$ such that the exceptional divisor $E$ supports a $\pi$-ample divisor and the support of $\mathrm{Supp}(E+\pi^{-1}_* \Delta)$ is a simple normal crossing divisor.
\end{definition}

Given a log resolution of $(X, \Delta)$ as above, we can write
\begin{equation}\label{log.discr.formula}
K_{X'}+\pi^{-1}_*\Delta+\sum b_iE_i=\pi^*(K_X + \Delta),
\end{equation}
where the $E_i$ are the irreducible components of $E$.

\begin{definition}
The \textit{log discrepancy} of $E_i$ with respect to $(X, \Delta)$ is $a(E_{i};X,\Delta) \defn 1-b_{i}$.
\end{definition}

Given a sub-log pair $(X, \Delta)$ and a geometric valuation $V$, we say that the valuation is exceptional if $V$ is not associated to any divisor on $X$. 
In this case, it is possible to find a log resolution $\pi : X' \to X$ such that $V$ is realized on $X'$ as the valuation associated to an exceptional 
prime Cartier divisor $D \subset X'$ (cf. \cite[Lemma 2.45]{koll-mor-book}). 

\begin{definition}
The log discrepancy
of $V$ is $a(V;X,\Delta) \defn a(D;X,\Delta)$.
\end{definition}
It is easy to verify that the definition of log discrepancy does not depend on the choice of the log resolution.\\
The center of $V$ on $X$, denoted $c_X(V)$ or $c_X(D)$, is defined as $\pi(D)$.
This notion is independent of the choice of the log resolution, too.

\begin{definition} 
The discrepancy of a sub-log pair $(X, \Delta)$ is
\[
{\rm discrep}(X, \Delta) \defn \inf \{{\rm a}(V; X, \Delta)\,| V \, {\rm divisorial \, valuation, \, exceptional\, over}\, X\}.
\]
For $Z \subset X$ an integral subvariety and $\eta_Z$ its generic point, we define
\begin{eqnarray}
\nonumber {\rm a}(Z; X, \Delta) &=& \inf_{V, \,c_X(V) \subseteq Z} {\rm a}(V; X, \Delta)\\
\nonumber {\rm a}(\eta_Z; X, \Delta) &=& \inf_{V, \,c_X(V) = Z} {\rm a}(V; X, \Delta).
\end{eqnarray}
\end{definition}

The log discrepancy of a divisorial valuations is the central object in the study of singularities of pairs.
It is a well known fact (cf. \cite{koll-mor-book}), that 
\[
0 \leq {\rm discrep}(X, \Delta) \leq \dim_\complex X \, {\rm or} \, {\rm discrep}(X, \Delta) = -\infty.
\]

The Minimal Model Program mainly focuses on studying those pairs whose log discrepancy is non-negative.

\begin{definition}
\label{sings.def}
A sub-log pair $(X,\Delta)$ is sub-kawamata log terminal (in short, sub-klt) (respectively sub-log canonical (sub-lc); sub-divisorial log terminal (sub-dlt)) if ${\rm discrep}(X, \Delta) > 0$ and $\lfloor \Delta \rfloor \leq 0$ (resp. ${\rm discrep}(X, \Delta) \geq 0$; if there exists a log resolution $\pi\colon X' \to X$ such that all exceptional divisors have log discrepancy $<1$).
\end{definition}

If $(X, \Delta)$ is a log pair then we highlight this property by removing the prefix {\it sub-} from all the notions of singularities defined in \ref{sings.def}.

\subsection{The non-klt locus}\label{nonlc.locus.sec}

\begin{definition} \label{nonlc.locus.def}
Let $(X, \Delta)$ be a sub-log pair and $Z \subset X$ an integral subvariety. 
Then, $Z$ is a non-kawamata log terminal center (in short, a non-klt center) if 
${\rm a}(\eta_Z; X, \Delta) \leq 0$.\\
The non-kawamata log terminal locus (non-klt locus) of the pair $(X, \Delta)$, $\nklt(X, \Delta)$, 
is the union of all the non-klt centers of $X$,
\[
\nklt(X, \Delta) \defn \bigcup_{\{Z | {\rm a}(\eta_Z; X, \Delta) \leq 0\}} Z.
\]
\noindent The non-log canonical locus (non-lc locus) of the pair $(X, \Delta)$, $\nlc(X, \Delta)$, is
\[
\nlc(X, \Delta) \defn \{X \ni p \, {\rm closed \, point} \, | {\rm a}(p; X, \Delta) =-\infty\}.
\]
\noindent $Z$ is a log canonical center (lc center) if $a(\eta_Z; X, \Delta) = 0$ and for a generic point $p \in Z, \, {\rm a}(p; X, \Delta) \geq 0$, 
i.e. $Z \nsubseteq \nlc(X, \Delta)$.
\end{definition}

\begin{remark}
Given a subvariety $Z \subset X$ for which ${\rm a}(\eta_Z; X, \Delta) <0$, 
then for every point $p \in Z, \, {\rm a}(p; X, \Delta) = -\infty$, as it is easy to verify 
by passing to a log resolution. Hence, the above definition 
of $\nlc(X, \Delta)$ is equivalent to the following alternative definition
\[
\nlc(X, \Delta) \defn \bigcup_{\{Z \subset X \, | \, {\rm a}(\eta_Z; X, \Delta)< 0\}} Z.
\]
\end{remark}

If we pass to a log resolution of $(X, \Delta)$, $\pi: X' \to X$ and write as in (\ref{log.discr.formula})
\[
K_{X'} + \Delta_{X'} = K_{X'}+\pi^{-1}_*\Delta+\sum b_iE_i=
\pi^*(K_X + \Delta) =  K_{X'}+ \sum_i b_i \Delta'_i,
\]
\noindent
then $\nklt(X, \Delta) = \pi({\rm Supp}(\sum_{b_i \geq 1} \Delta'_i))$ and
$\nlc(X, \Delta) = \pi({\rm Supp}(\sum_{b_i > 1} \Delta'_i))$.

The complement in $X$ of $\nklt(X, \Delta)$ is the biggest open set on which $\Delta$ has just sub-klt singularities and, analogously, the complement of 
$\nlc(X, \Delta)$ is the biggest open set of $X$ on which $\Delta$ has sub-lc singularities.

The divisor $\Delta_{X'}^{=1}$ is the source of lc centers of $\Delta$. 
It is easy to see (cf. \cite[Lemma 2.29]{koll-mor-book}) that all valuations of log discrepancy $0$ with 
respect to $(X, \Delta)$, not contained in $\nlc(X, \Delta),$ are given either by the components of 
$\Delta_{X'}^{=1}$ or by 
blowing up the strata of $\Delta_{X'}^{=1}$ and repeating the same procedure.
Hence, the lc centers are nothing but the closures of the lc centers for the pair 
$(X\setminus\nklt(X, \Delta), \Delta|_{X\setminus\nklt(X, \Delta)})$.

The union of the lc centers of $(X, \Delta)$ is a subvariety of $X$ (or a subscheme), but it carries 
a richer structure. It is in fact a subvariety stratified by the lc centers and it will be important for us 
to keep track of the strata. 

\begin{definition}
Let $(X, \Delta)$ be a log pair.
Given a lc center $W$ for $(X, \Delta)$, the total space of the stratification associated 
to $(X, \Delta)$ on W is given by
\[
\mathrm{Strat}(W, \Delta) \defn \bigcup_{\substack{W' \subsetneq W\\ W' \; {\rm lc \; center}}} W',
\]

the union of the log canonical centers contained in $W$.

\end{definition}

An important result about the structure of the non-klt locus, that we will need in the next sections 
of the paper, is the following connectedness theorem for negative maps, originally due to Shokurov.

\begin{theorem}\cite[Theorem 17.4]{koll-etal-book}\label{connect.thm}
Let $(X, \Delta)$ be a lc pair and let $\phi: X \to Y$ be a contraction of quasi-projective varieties. 
Assume that $-(K_X + \Delta)$ is $\pi$-nef and $\pi$-big. 
Then, every fiber of $\pi$ has a neighborhood (in the classical topology) in which the $\nklt(X, \Delta)$ is connected.
\end{theorem}

An extension of this theorem to the case of an lc pair $(X, \Delta)$ together with a contraction of quasi-projective varieties $\phi: X \to Y$ such that $-(K_X + \Delta)$ $\pi$-nef (but not necessarily $\pi$-big) -- and more generally in the context of generalized pairs -- has been announced in \cite{mio.connectedness}.

\section{Dlt modifications}
\label{dlt.mod.sect.}

\subsection{Adjunction: the different}

 When dealing with a pair $(X, \Delta)$ that is not snc easy examples show that 
the adjunction formula might need the introduction of a correction term. That is, given a component 
$D$ of $\Delta$
of coefficient $1$, it could happen that in the adjunction formula
\[
(K_X +D)|_{D} \neq K_D.
\] 
For more details on this, see \cite[\S 16]{koll-etal-book}.

When $(X, \Delta)$ is dlt, it is possible to modify the theory and 
obtain something analogous to the classical 
adjunction setting, that furthermore behaves well when restricting to higher codimension lc centers.

\begin{theorem}\label{diff.defn}
Let $(X,\Delta)$ be a dlt pair and $W\subset X$ be a lc center.
There exists on $W$ a naturally defined $\real$-divisor ${\rm Diff}^\ast_W \Delta \geq 0$
such that
\[
(K_X+\Delta)|_W\sim_\ratl K_W+{\rm Diff}^\ast_W \Delta
\]
and the pair $(W, {\rm Diff}^\ast_W \Delta)$ has dlt singularities. Moreover,
the non-klt locus of $(W, {\rm Diff}^\ast_W \Delta)$  is equal to the union of the lc centers of $\Delta$
strictly contained in $W$.
\end{theorem}

The divisor ${\rm Diff}^\ast_W \Delta$ can be defined inductively 
starting as in \cite[Sec.16]{koll-etal-book} 
from the case in which $W=D$ is a divisor. Then
\[
(K_X+D+(\Delta - D))|_D \sim_\ratl K_D+{\rm Diff}^\ast_D \Delta.
\]
Working inductively, when $(X, \Delta)$ is dlt, ${\rm Diff}^\ast_W\Delta$ is constructed analogously whenever
$W$ is an lc center of $(X, \Delta)$;
that is simply because the generic point of $W$ is an irreducible component of a complete intersection of components of $\lfloor \Delta \rfloor$.

\begin{definition}\label{diff.defn.true}
The divisor ${\rm Diff}^\ast_W\Delta$ from Theorem \ref{diff.defn} is called
the different of $\Delta$ on $W$.
\end{definition}

In general it is hard to obtain analogous adjunction results for higher codimensional lc centers.
Nonetheless, for an lc center which is minimal with respect to inclusion we have an analogous form of adjunction, originally due to Kawamata.

\begin{theorem}[Kawamata subadjunction] \cite{kawam91, amb05, fuj-gong12}
\label{kaw.subadj.thm} 
Let $(X, \Delta)$ be a log canonical pair and $W$ a minimal lc center.
Then there exists an effective divisor $\Delta_W$ on W s.t. $(W, \Delta_W)$ is klt and
\[
(K_X + \Delta)|_{W} \sim_{\real} K_W + \Delta_{W}.
\]
\end{theorem}

We refer the reader to \S \ref{subadj.sect} for a more detailed discussion of adjunction along lc centers of codimension higher than one.

\subsection{Existence of special dlt modifications}

An important fact, that will be needed multiple times in the following sections is that, starting with an lc pair, there always exists a crepant resolution giving a dlt pair. 
While this type of result has now been known for the past decade, see, for example, \cite{MR2955764}, in the case of Theorem \ref{main.thm}, we need a more refined resolution results that allows us to fully control not just the discrepancies of the codimension one part of the exceptional locus of the resolution. 

\begin{theorem}\label{dlt.mod.thm}
Let $(X, \Delta= \sum_i b_i D_i)$ be a log pair, $0 < b_i \leq 1$. 
Then there exists a $\ratl$-factorial pair $(Y, \Delta_Y = \sum_i b_i \Delta_i \geq 0)$ 
and a birational map $\pi\colon Y \to X$ with the following properties:

\begin{enumerate}

\item[(i)] $K_{Y} + \Delta_Y = \pi^*(K_X+\Delta)$; 

\item[(ii)] the pair $(Y, \Delta_Y' \defn \sum_{b_i < 1} b_i \Delta_i + \sum_{b_i \geq 1} \Delta_i)$ is dlt;

\item[(iii)] every divisorial component of ${\rm Exc}(\pi)$ appears in $\Delta_Y'$ with coefficient 1;

\item[(iv)] $\pi^{-1}(\nklt(X, \Delta))=\nklt(Y, \Delta_Y)= \nklt(Y, \Delta_Y')$.

\end{enumerate}

\end{theorem}

\begin{proof}

For the proof of (i), (ii), (iii) one can refer to \cite[3.10]{koll-kov10}.
Let $\pi_Z\colon (Z, \Delta_Z) \to X$ be a modification satisfying these properties.
Then 

\begin{equation}\label{dlt.eqna}
\Delta_Z = \Delta_Z^{<1} + \Delta_Z^{\geq1} = \sum_{i | b_i <1} b_i D_i + \sum_{i | b_i \geq 1} b_i D_i
\end{equation}
and $(Z, \Delta_Z^{<1})$ is a klt pair.
Moreover, as $K_{Z} + \Delta_Z = \pi_Z^*(K_X+\Delta)$,
\begin{equation}\label{dlt.pf.b}
K_{Z} + \Delta_Z^{<1} \sim_{\pi_Z, \real} -\Delta_Z^{\geq 1}.
\end{equation}
\newline
Therefore, we can run a relative $(K_{Z} + \Delta_Z^{<1})$-MMP over $X$, as $(Z,  \Delta_Z^{<1})$ is klt by construction;
this run of the MMP must terminate with a relatively minimal model $\psi\colon (Z, \Delta_Z^{<1}) \dasharrow (Z', \Delta_{Z'}^{<1}\defn \psi_* \Delta_Z^{<1})$, as $K_{Z} + \Delta_Z^{<1}$ is big/$X$, since $Z \to X$ is birational;
on the model $Z'$ the following conditions hold true:

\begin{itemize}

\item[a)] $(Z', \Delta_{Z'}^{<1})$ is a $\ratl$-factorial, klt pair;

\item[b)] $K_{Z'} + \Delta_{Z'}^{<1} + \Delta_{Z'}^{\geq 1} = \pi_{Z'}^*(K_X+\Delta)$, 
where $\Delta_{Z'}^{\geq1} \defn \psi_* \Delta_Z^{\geq1}$ and $\pi_{Z'}\colon Z' \to X$ is 
the structural map;

\item[c)] $K_{Z'} + \Delta_{Z'}^{<1}$ is $\pi_{Z'}$-nef and by \eqref{dlt.pf.b} the same 
holds for $-\Delta_{Z'}^{\geq1}$.

\end{itemize}
Properties a) and b) imply that $\nklt(X', \Delta_{Z'}^{<1} + \Delta_{Z'}^{\geq1}) = {\rm Supp}(\Delta_{Z'}^{\geq1})$.
In fact, the inclusion 
$\nklt(Z', \Delta_{Z'}^{<1} + \Delta_{Z'}^{\geq1}) \supseteqq {\rm Supp}(\Delta_{Z'}^{\geq1})$ follows form Definition
\ref{nonlc.locus.def}. To prove the other inclusion, let $W$ be a non-klt center not contained in 
Supp($\Delta_{Z'}^{\geq1}$). There exists a log resolution 
$r\colon (S, \Delta_S) \to (Z', \Delta_{Z'}^{<1} + \Delta_{Z'}^{\geq1})$ and a component $F_1$
of $\Delta_S$ whose coefficient is $\geq 1$ and $c_{Z'}(F_1) = W$. 
As $c_{Z'}(F_1) \nsubseteqq \Delta_{Z'}^{\geq1}$, it follows that
$a(F_1; Z', \Delta_{Z'}^{<1}) \leq 0$ as well, 
which is impossible as $(Z', \Delta_{Z'}^{<1})$ is klt.

Finally, take another dlt modification 
\[
\psi\colon (Y, \Delta_Y) \to (Z', \Delta_{Z'}^{<1} + \sum_{F_i \subset {\rm Supp}(\Delta_{Z'}^{\geq1})} F_i)
\] 
with properties (1), (2), (3) from the statement of the theorem. 
The divisor $-\psi^*(\Delta_{Z'}^{\geq1})$ will be a $\pi$-nef divisor, where $\pi= \pi_{Z'} \circ \psi$.
The support of $-\psi^*(\Delta_{Z'}^{\geq 1})$ contains all and only those components of $\Delta_Y$ of coefficient $\geq 1$. 
By negativity lemma, \cite[Lemma 3.39]{koll-mor-book},  $\pi\colon Y \to X$ satisfies condition (4) of the theorem.
\end{proof}

\section{Adjunction for higher codimensional lc centers}
\label{subadj.sect}

\subsection{Crepant log structures}
\label{crep.sect}

\begin{definition}
\cite[Def. 4.28]{koll-book13}
Let $Z$ be a normal variety.
A crepant log structure (respectively, dlt crepant log structure) on $Z$ is the datum of a normal log canonical (resp. dlt) pair $(Y, \Delta_Y)$ together with a contraction $f \colon Y \to Z$ such that $K_Y+\Delta_Y \sim_{f, \mathbb{Q}} 0$.
\newline
An irreducible subvariety $W \subset Z$ is a log canonical center (in short, lc center) of a crepant log structure $f \colon  (Y, \Delta_Y) \to Z$ if it is the image of an lc center $W_Y \subset Y$ of $(Y, \Delta_Y)$.
\end{definition}

Crepant log structures are very useful when studying adjunction for lc centers of codimension strictly greater than one for lc pairs.

Given an lc pair $(X, \Delta)$, let us fix a dlt modification $\pi \colon X' \to X$ of $(X, \Delta)$, $K_{X'}+ \Delta'= \pi^\ast(K_X+ \Delta)$.
We also fix an lc center $W \subset X$ with respect to $(X, \Delta)$ and denote by $W^\nu$ the normalization of $W$.
We choose $S \subset X'$ to be an lc center for $(X', \Delta')$ that dominates $W$. 
We consider the Stein factorization 
\[
\pi_{\vert S}\colon  S\stackrel{\pi_S}{\longrightarrow} W_S\stackrel{{\rm spr}_{W^\nu}}{\longrightarrow} W^\nu.
\]
By adjunction, the morphism $\pi_S \colon S \to W_S$ is a dlt crepant log structure over $W_S$ with datum  $(S, \Delta_S)$, where $\Delta_S := {\rm Diff}^\ast_S \Delta_{X'}$.

The following theorem, due to Koll\'ar,  shows that the contraction $\pi_S \colon S \to W_S$ already contains
all the relevant information in terms of the stratification of the non-klt locus of $(X, \Delta)$.

\begin{theorem}\cite[Cor. 4.42]{koll-book13}\label{dlt.crep.strat}
Let $\pi\colon(X',\Delta')\to Z$ be a dlt crepant log structure and $S \subset X'$ be an lc center, with $\pi(S)=W$.
We denote by $W^\nu$ the normalization of $W$.
We consider the Stein factorization 
\begin{equation*}
\pi\vert_{S}\colon  S\stackrel{\pi_S}{\longrightarrow} W_S\stackrel{{\rm spr}_W}{\longrightarrow} W
\end{equation*}
and let $\Delta_S\colon={\rm Diff}^\ast_S\Delta'$ be the different of $\Delta'$ on $S$.  
Then:
\begin{enumerate}
\item $\pi_S\colon (S, \Delta_S) \to W_S$ is a dlt, crepant log structure;

\item Given an lc center $Z_S\subset W_S$ for $\pi_S$, ${\rm spr}_W(Z_W)\subset W$ is an lc center for $\pi\colon(X',\Delta')\to Z$.
Every minimal lc center of $(S, \Delta_S)$ dominating $Z_S$ is also a minimal lc center of $(X',\Delta')$ that dominates $\pi(Z_W)$.

\item  For $Z\subset W$ an lc center of $\pi\vert_{S}\colon (S,\Delta_S)\to W$, every irreducible component of ${\rm spr}_W^{-1}(Z)$ is an lc center of $\pi_S\colon (S, \Delta_S) \to W_S$.
\end{enumerate}
\end{theorem}

With the notation adopted in Theorem \ref{dlt.crep.strat}, we will denote the total space of the stratification, induced on the lc center $W_S$ by the lc centers of $(W_S, \Delta_S)$ with
\[
\mathrm{Strat}(W_S, \Delta_S) \defn \bigcup_{\substack{W' \subsetneqq W,\\ W' \; {\rm lc 
\; center}}} \hspace{0.2cm} \bigcup_{\substack{V \, {\rm irreducible} \\
{\rm component \; of} \\ {\rm spr}^{-1}(W)}} V.
\]

\subsection{Canonical bundle formula}\label{cbf.subsect}

\begin{definition}\label{lc.triv.defn}\cite{fuj-gong.1210}
An lc-trivial fibration is the datum of a contraction
of normal varieties $\pi \colon Y \rightarrow Z$ and a pair $(Y,\Delta_Y)$ s.t.

\begin{enumerate}

\item[(i)] $(Y,\Delta_Y)$ has sub-lc singularities over the generic
point of $Y$, i.e., $\nlc(\Delta_Y)$ does not dominate $Z$ and $\Delta_Y$ could 
possibly contain components of negative coefficient;

\item[(ii)] 
$\mathrm{rank} \; \widetilde{\pi}_{\ast} \strutt_{\widetilde{Y}}(\lceil \mathbb{A}^{\ast}(Y,\Delta)\rceil) = 1$,
where $\tilde{\pi}=\pi \circ l$ and $l \colon \widetilde{Y} \to Y$ is a log resolution of $(Y,\Delta_Y)$.
$\mathbb{A}^{\ast}(Y,\Delta)$ is the b-divisor
whose trace on $\tilde{Y}$ is defined by the following equality
\[
K_{\widetilde{Y}} = \widetilde{\pi}^\ast(K_Y+\Delta_Y) + \sum_{a_i \leq -1} a_i D_i + \mathbb{A}^{\ast}(Y,\Delta)_{\widetilde{Y}}.
\]

\item[(iii)] there exist $r \in \natl$, a rational function $\phi \in k(Y)$
and a $\ratl$-Cartier divisor $D$ on $Y$ s.t.
\begin{eqnarray}\label{triv.overZ.eqn}
K_Y + \Delta_Y +\frac{1}{r}(\phi) = \pi^{\ast}D, & {\rm i.e.} \; \; K_Y +\Delta_Y \sim_{\pi, \ratl} 0. 
\end{eqnarray}

\end{enumerate}

 At times, we will denote an lc-trivial structure by $\pi\colon(Y, \Delta_Y) \to Z$. 

\end{definition}

\begin{remark}
A sufficient condition for (2) in definition \ref{lc.triv.defn} to hold is that $\Delta_Y$ is log canonical, in which case,
\[
\lceil \mathbb{A}^{\ast}(Y,\Delta_Y)_{\widetilde{Y}}\rceil=
\lceil K_{\widetilde{Y}}-\pi^{\ast}(K_Y+\Delta_Y)+\sum_{a(E,Y,\Delta_Y)=1} E\rceil
\]
is always exceptional over $Y$. Under this hypothesis, an lc-trivial fibration is also a crepant log structure.
\end{remark}

When working with lc-trivial fibrations, we are often interested in studying the relation between lc centers of the pair $(Y, \Delta_Y)$ and their images on $Z$, e.g., when working in the context of adjunction theory introduced in the previous section.

\begin{definition}
An integral subvariety $W\subset Z$ is an lc center of an lc-trivial fibration $\pi\colon Y\to Z$, if it is the image of an lc center $W_Y\subset Y$ for $(Y,\Delta_Y)$.
\end{definition}

We now introduce the main definitions and results related to the canonical bundle formula, that will be an essential tool in the proof of Theorem \ref{main.thm}.

\begin{definition}
Given an lc-trivial fibration $\pi\colon (Y, \Delta_Y) \to Z$ as above, 
let $T\subseteq Z$ be a prime divisor in $Z$.
The log canonical threshold of $\pi^{\ast}(T)$ with respect to the pair $(X,\Delta)$ is
\[
{\rm a}_T=\sup\{t\in\real | (Y,\Delta_Y+t\pi^{\ast}(T)) \mathrm{\;is \;lc\; over\;} T\}.
\]
We define the discriminant of $\pi \colon (Y,\Delta_Y) \to Z$ to be the divisor
\begin{equation*}
B_Z \defn \sum_{T}(1-{\rm a}_T)T.
\end{equation*}

\end{definition}

 It is easy to verify that the above sum is finite: 
a necessary condition for a prime divisor to have non-zero coefficient is to 
be dominated by some component of $B_Z$ of non-zero coefficient. 
There finitely many such components on $Y$. Hence, $B_Z$ is an $\mathbb{R}$-Weil divisor.

\begin{definition}
Let $\pi \colon (Y, \Delta_Y) \to Z$ be an lc-trivial fibration.
With the same notation as in equation (\ref{triv.overZ.eqn}),
fix $\phi\in k(Y)$ for which $K_Y + \Delta_Y +\frac{1}{r}(\phi) = \pi^{\ast}D$.
Then there is a unique divisor $M_Z$ for which the following equality holds
\begin{eqnarray}\label{cbf}
K_Y + \Delta_Y +\frac{1}{r}(\phi) &=& \pi^{\ast}(K_Z+B_Z+M_Z).
\end{eqnarray}

 The $\mathbb{Q}$-Weil divisor $M_Z$ is called the moduli part.

\end{definition}

When dealing with an lc-trivial fibration $\pi\colon(Y, \Delta_Y) \to Z$, we can pass to a higher birational model $Z'$ of $Z$ and take a resolution $Y'$ of the normalization of the main component of the fibre product $Y \times_Z Z'$ and form the corresponding cartesian diagram,
\begin{equation*}
\xymatrix{
Y \ar[d]_{\pi} & Y'
\ar[l]_{r_{Y}} \ar[d]^{\pi'}\\
Z & Z' \ar[l]_{r}.}
\end{equation*}

 By base change, we get a new pair, $(Y', \Delta_{Y'})$, from the formula 
\[
K_{Y'} +\Delta_{Y'} = r_{Y'}^\ast(K_Y+\Delta_Y).
\]

It follows from the definition that, under these hypotheses,  $\pi' \colon (Y', \Delta_{Y'}) \to Z'$
will be an lc-trivial fibration as well, allowing to compute $B_{Z'}$ and $M_{Z'}$.

 The discriminant and the moduli divisor have a
birational nature: they are b-divisors, as their definition immediately implies that
\[
r_\ast B_{Z'}= B_Z, \; {\rm and \;} \; r_\ast M_{Z'}=M_Z. 
\]
As they are b-divisors, we will denote them using the symbols $\mathbb{B}$ and $\mathbb{M}$, respectively.

 Fujino and Gongyo proved, generalizing results of Ambro,
 that these divisors have interesting features.
 
\begin{theorem}\label{nefness} (\cite{fuj-gong.1210}, \cite{amb05})
Let $\pi \colon (Y,\Delta_Y)\rightarrow Z$ be an lc-trivial fibration.
There exists a birational model
$Z'$ of $Z$ on which the following properties are satisfied:

\begin{itemize}

\item[(i)] $K_{Z'}+\mathbb{B}_{Z'}$ is $\mathbb{Q}$-Cartier, 
and $\mu^{\ast}(K_{Z'}+\mathbb{B}_{Z'}) = K_{Z''}+\mathbb{B}_{Z''}$
for every higher model $\mu\colon Z''\rightarrow Z'$. 

\item[(ii)] $\mathbb{M}_{Z'}$ is nef and $\mathbb{Q}$-Cartier. 
Moreover, $\mu^{\ast}(\mathbb{M}_{Z'}) = \mathbb{M}_{Z''}$ for every
higher model $\mu\colon Z''\rightarrow Z'$. More precisely, it is b-nef and good, i.e., 
there is a contraction $h \colon Z' \to T$ and $\mathbb{M}_{Z'} = h^\ast H$, for some $H$ big and nef on 
$Z'$.

\end{itemize}
\end{theorem}

When the model $Z'$ satisfies both conditions in the theorem, we say that $\mathbb{B}$ and $\mathbb{M}$ descend to $Z'$.

\subsection{Springs and sources}
\label{adj.high.codim.sect}

One of the main reasons to study crepant log structures and lc-trivial fibrations comes from resolutions and adjunction.
Let $(X, \Delta)$ be an lc pair and $W \subset X$ an lc center.
In the purely lc case, when $(X, \Delta)$ is not dlt, the adjunction theory along $\nklt(X, \Delta)$ is not as easily determined as in Theorem \ref{diff.defn}. 
Nonetheless, Theorem \ref{dlt.mod.thm} shows that it is always possible to pass to a dlt pair crepant to the original one. 
Let $\pi\colon X' \to X$ be a dlt modification as in the Theorem \ref{dlt.mod.thm}, with $K_{X'}+\Delta'= \pi^\ast(K_X+\Delta_X)$.

Let $S$ be a log canonical center of $\Delta_{X'}$, i.e., an irreducible component of intersections 
of components of coefficient $1$ such that $W$ is the image of $S$ on $X$. 
As we saw in Theorem \ref{dlt.crep.strat}, the datum of the Stein factorization of the map $S \to W$, 
\begin{equation}
\label{stein.fact.spring}
\pi\vert_{S}\colon  S\stackrel{\pi_S}{\longrightarrow} W_S\stackrel{{\rm spr}_W}{\longrightarrow} W
\end{equation}
yields a dlt crepant log structure which is also an lc-trivial fibration. 
If the lc center $S$ is chosen to be minimal among those dominating $W$, then the singularities of $(S, {\rm Diff}^\ast_S \Delta')$ are actually of klt type over the generic point of $W$.

Under this minimality condition on the lc center $S$, Koll\'ar, cf. \cite[\S 4.5]{koll-book13}, proved that the isomorphism class of the variety $W_S$ over $W$ in \eqref{stein.fact.spring} is independent of the choice of $S$. 
Moreover, he also showed that for any two pairs
$(S_1, \Delta_{S_1})$, $(S_2, \Delta_{S_2})$ of minimal lc centers dominating $W$ the pairs $(S_i, \Delta_{S_i})$ are crepant birational, that is, $S_1$ and $S_2$ are birational and there exists a common resolution $p_i\colon T \to S_i, $, $i=1, 2$ such that
\[
p_1^\ast(K_{S_1}+\Delta_{S_1})=p_2^\ast(K_{S_2}+\Delta_{S_2}).
\]

This prompts the following definition.
\begin{definition}\cite[Thm./Def. 4.45]{koll-book13}
With the notation of this section, let $S$ be
an lc center of $(X', \Delta')$ minimal with respect to inclusion among the lc centers $T$ with $\pi(T)=W$. 
We call the pair $(S, \Delta_S={\rm Diff}^\ast_S\Delta_Y)$ a source of $W$.
\newline
The normal variety $W_S$ appearing in the Stein factorization of the morphism $\pi|_S\colon S \to W$, \eqref{stein.fact.spring}, is called the spring of $W$.
\end{definition}

\section{Mori hyperbolicity}\label{mh.sect}

The following is a generalization of the definition of Mori hyperbolicity that originally appeared in \cite{lu-zh}.

\begin{definition}\label{m.hyp.def}
Let $(X, \Delta= \sum_{i}b_i D_i), \; 0< b_i \leq 1$ be a log pair. We say that $(X, \Delta)$ is a Mori 
hyperbolic pair if 
\begin{enumerate}
\item there is no non-constant morphism $f\colon \mathbb{A}^1 \to X \setminus \nklt(X, \Delta)$; 
\item for any $W \subset X$ lc center, there is no non-constant morphism 
\[
f\colon \mathbb{A}^1 \to W \setminus \{(W \cap \nlc(\Delta)) \cup {\rm Strat}(W, \Delta)\}.
\]
\end{enumerate}

\end{definition}

The following inductive result is already implicitly contained in \cite[\S 4]{lu-zh}: it will be the starting point of our approach to the proof of Theorem \ref{main.thm}. 
We restate it here for the reader's convenience since it does not appear in \cite{lu-zh} in this generality. 

\begin{proposition}\label{mh.nonlc.dlt.prop}
Let $(X, \Delta = \sum_i b_i D_i \geq 0)$ be a normal, projective, $\ratl$-factorial log pair such that $(X, \Delta' =\sum_{i | b_i < 1} b_i D_i + \sum_{i | b_i \geq 1} D_i)$ is 
dlt. 
\newline
Suppose that $K_X + \Delta$ is nef when restricted to  ${\rm Supp}(\Delta^{\geq 1})$.
Then,
\begin{enumerate}
\item[(i)] either $K_X+ \Delta$ is nef or 
\item[(ii)] there exists a non-constant morphism $f\colon \mathbb{A}^1 \to (X \setminus \nklt(X, \Delta))$.
\end{enumerate}
\end{proposition}

\begin{proof}
Suppose $K_X + \Delta$ is not nef. 
Then there exists a $(K_X+ \Delta)$-negative extremal ray $R$ in the cone of effective curves $\NE(X)$. 
Since $K_X+\Delta$ is nef along $\nklt(X, \Delta)$, $R$ is both a $(K_X+ \Delta')$-negative and a $(K_X+ \Delta^{<1})$-negative extremal ray, as $\nklt(X, \Delta)={\rm Supp}(\Delta^{\geq 1})$, by the assumptions of the proposition.
In particular, there exists an extremal contraction $\mu\colon X \to S$ associated to $R$.  
\newline
As $R$ does not contain classes of curves laying in $\nklt(X, \Delta)$, $\mu$ induces a finite morphism 
when restricted to $\nklt(X, \Delta)$. 
Thus, the $\ratl$-factoriality of $X$ implies that we are in either of these three cases:
\begin{enumerate}
\item[1)] $\mu$ is a Mori fibre space and all the fibres are one dimensional;

\item[2)] $\mu$ is birational and the exceptional locus does not intersect $\nklt(X, \Delta)$;

\item[3)] $\mu$ is birational and the exceptional locus intersects $\nklt(X, \Delta)$.

\end{enumerate}
\noindent 
As $\mu$ is a $(K_X+ \Delta^{<1})$-negative fibration and $K_X+ \Delta^{<1}$ is klt, then all of its fibres are rational chain connected, by \cite[Corollary 1.5]{hac-mck05}. 
Moreover,
\begin{eqnarray}\label{vanishing.OW.eqn}
R^1\mu_* \strutt_W = 0,
\end{eqnarray} 
\noindent
by relative Kawamata-Viehweg vanishing \cite[page 150]{laz2.book}. 
Thus, Theorem \ref{connect.thm} implies that $\nklt(X, \Delta') = \nklt(X, \Delta)$ is connected in a neighborhood of every fibre.
\newline
In case 1), the generic fibre of $\mu$ is a smooth projective rational curve. Theorem \ref{connect.thm} implies that the generic fibre intersects $\nklt(X, \Delta)$  in at most one point. 
This concludes the proof in case 1).
\newline
In case 2), as the fibres of $\mu$ are rationally chain connected, there exists a rational projective curve contained in $X \setminus \nklt(X, \Delta)$. 
This concludes the proof in case 2).
\newline
In case 3), we claim that the positive dimensional fibres are chains of rational curves.
To prove this claim, it suffices to show that each positive dimensional fiber is one-dimensional; rationality then follows as above from the results of \cite{hac-mck05}.
Assume that there is a positive dimensional fibre $F$ of dimension $>1$. 
By \cite[Corollary 1.5]{hac-mck05}, $F$ is covered by $(K_X+\Delta)$-negative rational curves. 
Hence, $F$ must intersect $\nklt(X, \Delta)$ or, else, we are in case $(ii)$ of the statement of the proposition.
As $X$ is $\mathbb{Q}$-factorial and $\dim F > 1$, then $\dim (F \cap \nklt(X, \Delta)) = \dim F -1 \geq 1$.
Hence, $\nklt(X, \Delta)$ contains curves that are contracted by $\mu$, but such curves are $(K_X+\Delta)$-negative since $\rho(X/S)=1$, contradicting the nefness of $K_X+\Delta$ along $\nklt(X, \Delta)$.
Therefore, from this claim it follows that, by the vanishing in \eqref{vanishing.OW.eqn} above, the positive dimensional fibers are trees of smooth rational curves. 
By Theorem \ref{connect.thm}, $\nklt(X, \Delta)$ intersects this chain in at most one point. 
In particular, there exists a complete rational curve $C$ such that $C \cap (X \setminus \nklt(X, \cdel)) = f(\mathbb{A}^1)$, where $f$ is a non-constant morphism.
This concludes the proof in case 3).
\end{proof}

In the case of a general log pair, using dlt modifications we get the following criterion, which will be fundamental in the proof of Theorem \ref{main.thm}.

\begin{corollary}\label{fund.cor}
Let $(X, \Delta= \sum_i b_i D_i \geq 0), 0 <b_i \leq 1$ be a log pair. 
Assume that there is no non-constant morphism $f\colon \mathbb{A}^1 \to X \setminus \nklt(X, \Delta)$.
\newline
Then $K_X + \Delta$ is nef if and only if $K_X + \Delta$ is nef when restricted to $\nklt(X, \Delta)$.
\end{corollary} 

\begin{proof}
Nefness of $K_X + \Delta$ immediately implies nefness of its restriction to every subscheme of $X$.
Hence, we just have to prove the converse implication.
\\
Let $\pi\colon (X', \Delta_{X'}) \to (X, \Delta)$ be a dlt modification for $(X, \Delta)$ as in Theorem  \ref{dlt.mod.thm}. We can reduce to proving nefness for $K_{X'}+\Delta_{X'}$. 
As $\pi(\nklt(X', \Delta_{X'})) = \nklt(X, \Delta)$, $K_{X'} + \Delta_{X'}$ is nef when restricted to $\nklt(X', \Delta_{X'})$.
\\
Suppose $K_{X'} + \Delta'$ is not nef. By the proposition,
there exists a non-constant morphism $f\colon \mathbb{A}^1 \to (X' \setminus \nklt(X', \Delta'))$. 
This contradicts the assumption in the statement of 
the corollary, as the properties of dlt modifications imply that the image of $\pi \circ f$ lies in 
$X\setminus \nklt(X, \Delta)$.
\end{proof}

Let us notice that in the above corollary, we did not impose any condition on the singularities 
of $\Delta$, besides the coefficients being in $[0,1]$.

\section{Proof of theorem \ref{main.thm}}\label{proof.sect}

 We will work inductively on the strata of $\nklt(X, \Delta)$. 
Namely, we will prove that $K_X + \Delta$ is nef when restricted to every stratum of $\nklt(X, \Delta)$.
As the union of all the strata is the non-klt locus itself, 
the theorem will follow from Corollary \ref{fund.cor}.\vspace{0.2cm}

 \noindent \textbf{\textit{Step 1. Start of the induction: the case of minimal lc centers.}}

 When $W$ is a minimal lc center with respect to inclusion, then nefness of $(K_X+\Delta)|_{W}$ follows from Theorem \ref{kaw.subadj.thm}.
Recall that by that result there exists an effective divisor $\Delta_W$ on $W$ such that $(W, \Delta_W)$ is klt and $(K_X + \Delta)|_{ W} \sim_{\real} K_W + \Delta_{W}$.
By definition of Mori hyperbolicity $W$ does not contain rational curves; hence, it follows that $K_W + \Delta_{W}$ must be nef by the Cone theorem.
\vspace{0.1cm}

\noindent \textbf{\textit{Step 2. Moving the computation to the spring of W.}}

 We assume now that $W$ is no longer minimal and  that $K_X+\Delta$ is
nef when restricted to any other stratum $W'$ strictly contained in $W$.
 Recall the following notation
\[
\mathrm{Strat}(W, \Delta) = \bigcup_{\substack{W' \subsetneqq W,\\ 
W' \; {\rm lc \; center}}} W'
\]
to indicate the union of all substrata contained in $W$.

Let us fix a dlt modification of $(X, \Delta), \; \pi\colon (X', \Delta') \to (X, \Delta)$.
We also fix a non-klt center $W \subset X$ and let $S \subset X'$ be an lc center, minimal among
those dominating $W$. 
Let us consider the Stein factorization 
\[
\pi_{\vert S}\colon  S\stackrel{\pi_S}{\longrightarrow} W_S\stackrel{{\rm spr}_W}{\longrightarrow} W.
\]
The variety $W_S$ is normal, projective and it is naturally equipped with the $\real$-divisor 
\[
L \colon={\rm spr}_W^*(K_X+\Delta).
\]
The morphism $\pi_S\colon S \to W_S$ is an lc-trivial fibration with respect to $\Delta_S = {\rm Diff}^\ast_S \Delta_{X'}$ on $S$, as we saw in \S \ref{adj.high.codim.sect}; it is a dlt log crepant structure, too. 
Proving nefness of $(K_X+\Delta)_{\vert W}$ is equivalent to proving nefness of $L$; 
therefore, we can assume that $L$ is nef on $\mathrm{Strat}(W_S, \Delta_S)$ since
\[
\mathrm{Strat}(W_S, \Delta_S)= {\rm spr_W}^{-1}({\rm Strat}(W, \Delta)),
\]
by (3) in Theorem \ref{dlt.crep.strat}.
\newline
Hence, without loss of generality, we could substitute the triple $(W, (K_X+\Delta)|_{ W}, \mathrm{Strat}(W, \Delta))$ with the triple $(W_S, L, \mathrm{Strat}(W_S, \Delta_S))$.
In fact, if $L$ is not nef, then we will show that there exists a non-constant morphism $f\colon \mathbb{A}^1 \to W_S \setminus \mathrm{Strat}(W_S, \Delta_S)$. 
By Theorem \ref{dlt.crep.strat}, it follows that there exists a non-constant morphism $f'\colon \mathbb{A}^1 \to W \setminus \mathrm{Strat}(W, \Delta)$, violating the Mori hyperbolicity assumption for $W$.\\
To simplify the notation, we will denote $W_S$ by $W$, and $\mathrm{Strat}(W_S, \Delta_S)$ by ${\rm Strat}(W, \Delta_S)$.
\vspace{0.2cm}

\noindent \textbf{\textit{Step 3. Constructing a good approximation for $L$ on $W$.}}

By the results of Section \ref{subadj.sect}, there exist sufficiently high birational models $S'$ of $S$ and
$W'$ of $W$ together with a commutative diagram
\begin{eqnarray}\label{fbr.prod.diagr.2}
\xymatrix{
S \ar[d]_{\pi_S} & S' \ar[l]_{r_{S'}} \ar[d]^{\pi_{S'}}\\
W & W' \ar[l]_{r}}
\end{eqnarray}
having the following properties:

\begin{enumerate}

\item $r^*(L) = K_{W'} + \mathbb{B}_{W'} + \mathbb{M}_{W'}$;

\item $(W', \mathbb{B}_{W'})$ is snc and sublc, i.e., $\mathbb{B}_{W'}$ is not necessarily effective;
 
\item $K_{W'} + \mathbb{B}_{W'}$ descends to $W'$ and $\mathbb{M}_{W'}$ is nef and abundant.

\item $(S', \Delta_{S'})$ is a sublc pair, where $K_S'+\Delta_{S'}=r_{S'}^\ast(K_S+\Delta_{S})$;

\end{enumerate}
In this context, we compare singularities of $(W', \mathbb{B}_{W'}$) with those 
of the original pair $(W, \Delta)$. 
\begin{lemma}\label{same.nklt.loc.lemma}
With the above notation and hypotheses, we have that $r(\nklt(W', \mathbb{B}_{W'})) 
= {\rm Strat}(W, \Delta_S)$.

\end{lemma} 

\begin{proof}
We know that 
$r_{S'}(\nklt(X', \Delta_{S'})) = \nklt (S, \Delta_S)$ and $\pi_S(\nklt(S, \Delta_S)) = {\rm Strat}(W, \Delta)$. 
As the diagram in (\ref{fbr.prod.diagr.2}) commutes, we need to prove that 
$\pi_{S'}(\nklt(S', \Delta_{S'})) = \nklt(W', \mathbb{B}_{W'})$.
The definition of $\mathbb{B}_{W'}$ implies that every stratum of 
$ \nklt(W', \mathbb{B}_{W'})\subset W'$ is dominated by a stratum of $\nklt(S', \Delta_{S'})$, hence 
$\nklt(W', \mathbb{B}_{W'}) \subset \pi_{S'}(\nklt(S', \Delta_{S'}))$. 
The opposite inclusion is also true, as given a stratum of $\nklt(S', \Delta_{S'})$, up to going to 
higher models of $W'$ and $S'$, we can suppose that $D$ is a divisor whose image $D'$ on $W'$ is a divisor, too. In this case, by the definition of $\mathbb{B}_{W'}$  and since it descends to $W'$, $D' \subset \nklt(W', \mathbb{B}_{W'})$. 
Thus, $\nklt(W', \mathbb{B}_{W'}) \supset \pi_{S'}(\nklt(S', \Delta_{S'}))$.
\end{proof}
\noindent As proving that $L$ is nef is equivalent to proving that, for any given ample Cartier divisor 
$A$ on $W$ and any given $\epsilon>0$, 
$L + \epsilon A$ is nef, we focus on the divisor
\begin{eqnarray}\label{eqn.cbf.2}
r^*(L+ \epsilon A)
= K_{W'} + \mathbb{B}_{W'} + \mathbb{M}_{W'}+ r^*(\epsilon A).
\end{eqnarray}
By construction, we can assume that there exists an effective divisor $E$ supported 
on the exceptional locus of $r$ and $-E$ is relatively ample over $W$.
Hence, there exists a positive number $\theta_\epsilon\ll \epsilon$, such that
for any $0<\delta \leq \theta_\epsilon$, $\mathbb{M}_{W'}+ r^*(\epsilon A) - \delta E$ is an 
ample divisor on $W'$.

\begin{lemma}\label{step4.lemma}
For every $\epsilon > 0$, there is a suitable choice of  $\delta$ and of an
effective $\real$-divisor
$Q_\epsilon \sim_\real \mathbb{M}_{W'}+ r^*(\epsilon A) - \delta E$ 
for which the following equalities hold
\[
\mathrm{Nklt}(W', \mathbb{B}_{W'} + \delta E + Q_\epsilon) = 
\nklt(W', \mathbb{B}_{W'} + \delta E) = \nklt(W', \mathbb{B}_{W'}).
\]
With this notation, 
\begin{eqnarray}\label{eqn.cbf.4}
r^*(L+ \epsilon A) \sim_{\real} K_{W'} + \mathbb{B}_{W'} + \delta E + Q_\epsilon.
\end{eqnarray}
\end{lemma}

\begin{proof}
The first equality is a consequence of \cite[Proposition 9.2.26]{laz2.book}, once 
we choose $\delta$ small enough so that $Q_\epsilon$ is ample. 
The second equality follows immediately from the fact that we can choose $\delta$ 
to be arbitrarily small, since $(W', \mathbb{B}_{W'})$ is snc and sublc.
\end{proof}

\vspace{0.2cm}

\noindent \textbf{\textit{ Step 4. End of the proof.}}

\noindent Using Lemma \ref{step4.lemma}, we define a new divisor on $W$
\[
\Gamma_\epsilon \defn r_{\ast}(\mathbb{B}_{W'} + \delta E + Q_\epsilon).
\]
The pair $(W, \Gamma_\epsilon)$ is a log pair and its coefficients are real numbers in $[0,1]$.
By construction, those coefficients in $\mathbb{B}_{W'} + \delta_\epsilon E + Q_\epsilon$ that are
strictly larger than $1$ were those of components that are exceptional over $W$. Also,
$L+ \epsilon A \sim_\real K_W + \Gamma_\epsilon$ and we are reduced to
proving nefness for $K_W+\Gamma_\epsilon, \; {\rm for} \; \epsilon \ll 1$. The pair
$(W, \Gamma_\epsilon)$ fails to be lc but $\nklt(W, \Gamma_\epsilon) = 
\mathrm{Strat}(W, \Delta_S)$, by Lemma \ref{same.nklt.loc.lemma} and Lemma \ref{step4.lemma}.
Moreover, $K_W + \Gamma_\epsilon$ is nef, more precisely ample, when restricted to 
its non-klt locus. Hence, it is nef on $W$ by Corollary \ref{fund.cor}.
Since this holds for arbitrary choice of $\epsilon>0$, it follows that $L$ is nef on $W$,
terminating the proof of the inductive step and of the theorem.
\QEDB

\begin{remark}
In the course of the proof of \ref{main.thm}, we have shown the following (very) weak version of (quasi log canonical) subadjunction.
Surely, this is not the most desirable version of subadjunction that is expected to hold, as we explain below.

\begin{theorem}\label{nonlc.subadj.thm}
Let $(Y, \Delta)$ be a log canonical pair and $\pi\colon Y \to Z$ be an lc trivial fibration. 
Let $A$ be an ample divisor on $Z$.

Then for all $ \epsilon, \delta >0$, there exists an effective divisor $\Gamma_{\epsilon, \delta}$, 
with coefficients in $[0,1]$ satisfying the linear equivalence relation
\[
K_Z+\mathbb{B}_Z+\mathbb{M}_Z+\epsilon A \sim_\real K_Z+\Gamma_{\epsilon, \delta}.
\]

The pair $(Z, \Gamma_{\epsilon, \delta})$ is not log canonical, but there exists a log resolution 
$\pi\colon Z' \to Z$ such that the log discrepancy of the $\pi$-exceptional divisors is bounded 
below by $-\delta$, i.e.
\[
a(E; Z, \Gamma_{\epsilon, \delta})> -\delta, \; \textrm{for every} \; E \subset Z' \; \textrm{prime divisor exceptional over}\; Z.
\]
\end{theorem}
It was conjectured in \cite{1608.02997} and later proved in \cite{1807.04847} that an analogous of the canonical bundle formula holds in the category of generalized pairs.
The moduli b-divisor, $\mathbb{M}$, is expected to be semi-ample on a sufficiently high birational model of $Z$.
That would easily imply that, for a certain choice of $\mathbb{M}_Z$,  $(Z, \mathbb{B}_Z+\mathbb{M}_Z)$ is log canonical. 
If that were to be true, the proof of Theorem
\ref{main.thm} could be considerably simplified. In fact, $L$ would be linearly 
equivalent to the lc divisor $K_Z+\mathbb{B}_Z+\mathbb{M}_Z$ and 
\[
\nklt(Z, \mathbb{B}_Z+\mathbb{M}_Z)= \nklt(Z, \mathbb{B}_Z)= {\rm Strat}(W, \Delta).
\]
\end{remark}

In the proof of the Theorem \ref{main.thm} we showed that if $K_X+\Delta$ is not nef,
there is a non-constant morphism $f\colon \mathbb{A}^1 \to X$ whose image is contained in an lc center $W \subset X$ and it
does not intersect the lc centers strictly contained in $W$. In particular, from the inductive procedure used in the proof, we 
see that it is possible to select $W$ to be a minimal lc center among those on which the restriction of $K_X+\Delta$ is not nef. Hence, we obtain as a consequence we obtain the following generalized version of the Cone Theorem.

\begin{theorem}\label{generalized.cone.thm2}
Let $(X, \Delta)$ be a log canonical pair. 
Then there exist countably many $(K_X+\Delta)$-negative rational curves $C_i$ 
such that
\[
\NE(X)= \NE(X)_{K_X+\Delta \geq 0} + \sum_{i \in I} \real_{\geq 0}[C_i].
\]
Moreover, one of the two following conditions hold:
\begin{itemize}
\item $C_i \cap (X\setminus \nklt(X, \Delta))$ contains the image of a non-constant morphism \newline{$f\colon\mathbb{A}^1 \to X\setminus \nklt(X, \Delta)$};
\item there exists an lc center $W \subset X$ such that $C_i \cap (W\setminus \mathrm{Strat}(W, \Delta))$ contains the image of a non-constant morphism $f\colon \mathbb{A}^1 \to (W\setminus \mathrm{Strat}(W, \Delta))$.
\end{itemize}
\end{theorem}

 In an attempt to expand the above results to arbitrary singularities,
the following questions appear quite natural.

\begin{question}\label{gen.case.question}
Let $(X, \Delta = \sum b_i D_i \geq 0), 0 < b_i \leq 1$, be a Mori hyperbolic log pair. 
Assume $K_X+ \Delta$ is nef when restricted to $\nlc(X, \Delta)$. 
Is $K_X+ \Delta$ nef? Is it possible to drop the assumption $0 < b_i \leq 1$?
\end{question}

Most of the proof of Theorem \ref{main.thm} applies to the case of varieties with worse singularities than log canonical,  through the language and techniques of quasi log varieties introduced in \cite{amb03}. 
It seems that, in order to finish the proof, one would have to prove a stronger version of the Bend and Break Lemma. 
Assuming that the pair $(X, \Delta)$ has dlt support as in Theorem \ref{dlt.mod.thm}, one would need not just to be able to deform a $(K_X+\Delta)$-negative curve, but to do so while also keeping fixed the support of the intersection of the curve the support of the non-klt locus of $(X, \Delta)$.
In this way, by iterating this procedure, one should eventually arrive to create a copy of the affine line that is contained in the complement of the non-klt locus or in some strata of it.
Unfortunately, we are not able to prove such a result at this time, hence the above question remains still open. 
Some results in this direction were recently proved by McQuillan and Pacienza in \cite{mcq-pac}, for quotient singularities.

To address Question \ref{gen.case.question}, one could mimic the same proof as for 
Theorem \ref{main.thm}. Namely, starting with a log pair $(X, \Delta)$ such that the coefficients of 
$\Delta$ are in $[0, 1]$, no matter what the singularities of $\Delta$ are, it is sufficient to prove that
$K_X+\Delta$ is nef on $\nklt(X, \Delta)$, by Corollary \ref{fund.cor}. As there is very little 
control on the non lc locus of $\Delta$ (cf. \cite[Theorem 0.2]{amb03}), it seems 
inevitable to assume the nefness for the restriction of $K_X+\Delta$. In this setting, 
the formalism of the canonical bundle formula is not available anymore, but in order to study adjunction or
just the restriction of $K_X+\Delta$ to lc centers of $\Delta$, the formalism of log varieties 
can be used (cf. \cite{amb03} and \cite{fuj09}). Again, working by induction, one can restrict to 
a given stratum, $W$, and assume that nefness is known for the smaller strata and the intersection with the 
non-lc locus. Assuming by contradiction that $(K_X+\Delta)|_{W}$ is not nef, then we can find a contraction
morphism $\pi\colon W \to S$ which contracts curves with $(K_X+\Delta)$-negative class in a given extremal ray
contained in $\NE(X)$. It is not hard to prove that the fibres of $\pi$ will contain rational curves.
The hard part is to prove that it is possible to deform one of these curves to a rational curve whose 
normalization supports the pull-back of $\Delta$ at most one point. The classical tool to deform curves is 
surely the Bend and Break Lemma, although in this case, we need not only to be able to deform a curve, 
but also we would like to be able to control its intersection with the components of $\Delta$. Hence, ideally,
one would like to prove a stronger version of the Bend and Break Lemma that makes the above construction
possible.

\section{Ampleness and pseudoeffectiveness for Mori hyperbolic pairs}\label{ample.sect}
When dealing with Mori hyperbolic pairs, in the dlt case, 
one can actually go further and give criteria for the ampleness of $K_X+\Delta$ as described in 
Theorem \ref{ample.thm} in the Introduction. Such criteria are modeled along the lines of 
the classical Nakai-Moishezon-Kleiman criterion which we
recall here in the version for $\mathbb{R}$-divisors due to Campana and Peternell.

\begin{theorem}[Campana-Peternell] \cite[Thm. 2.3.18]{laz1.book}
\label{kleiman.positivity.thm} 
Let $X$ be a proper variety and let $D$ be an 
$\mathbb{R}$-Cartier $\mathbb{R}$-divisor on $X$.

Then $D$ is ample on $X$ if and only if for every proper subvariety variety $Y \subseteqq X$
\[
\int_Y D^{\dim Y}> 0.
\]
\end{theorem}
We will also need the following definition.
\begin{definition}
Let $(X, \Delta)$ a log canonical pair. An $\real$-divisor $D$ is log big (with respect to $(X, \Delta)$) if $D$ is big
and $D|_{W}$ is big for any lc center $W$ of $\Delta$.
\end{definition}

\begin{proposition}\label{nakai.crit.m.hyp.thm}
Let $(X, \Delta)$ be a log canonical pair. Then the following are equivalent:
\begin{enumerate}

\item[(1)] the divisor $K_X+\Delta$ is ample;

\item [(2)] the divisor $K_X+\Delta$ is big, its restriction to $\nklt(X, \Delta)$ is ample and $K_X+\Delta$ has strictly positive degree on 
every rational curve intersecting $X \setminus \nklt(X, \cdel)$.

\end{enumerate}
If $(X, \cdel)$ is dlt, then the above conditions are also equivalent to:

\begin{enumerate}

\item[(3)] the divisor $K_X+\Delta$ is nef and log big and it has strictly positive degree on every rational curve.

\end{enumerate} 

\end{proposition}

\begin{remark}
The assumption on the bigness of $K_X+\Delta$ in the proposition is necessary as the following example shows.

Let $E$ be a curve of genus $0$. Then $K_E \sim 0$ and the pair $(E, 0)$ is terminal (hence, log canonical) with empty non-klt locus. The curve $E$ clearly does not contain rational curves, nonetheless $K_E$ is not ample.
\end{remark}
\begin{proof}[Proof of Proposition \ref{nakai.crit.m.hyp.thm}.]
Clearly condition (1) implies conditions (2) and (3). 

Condition (2) implies that $K_X+\Delta$ is nef. In fact, by the Cone Theorem, an extremal 
ray contained in $\NE(X)$ on which $K_X+\Delta$ is negative is spanned by the class of a rational 
curve $C \subset X$. As $K_X+\Delta$ is ample along $\nklt(\Delta)$, $C$ must intersect $X \setminus 
\nklt(X, \Delta)$, which gives a contradiction. \newline
Thus, $K_X+\cdel$ is big and nef and it is ample along $\nklt(X, \Delta)$. It follows that $K_X+\Delta$ is semiample, by \cite[Thm. 4.1]{fuj09}. The corresponding morphism is either an 
isomorphism or it has to contract some rational curves intersecting $X \setminus \nklt(X, \Delta)$ as 
implied by \cite[Thm. 1.2]{hac-mck05}. But this also gives a contradiction, as the intersection of $K_X+\cdel$
with such curves must be strictly positive. Then (2) implies (1).

 Let us prove that (3) implies (2). Since $K_X+\Delta$ is nef and log big, it is also semiample. By induction on the dimension and using Theorem \ref{diff.defn}, it follows
that $K_X+\cdel$ is ample along $\lfloor \cdel \rfloor$, which concludes the proof.
\end{proof}

\begin{theorem} \label{nakai.crit.m.hyp.cor}
Let $(X, \Delta)$ be a Mori hyperbolic log canonical pair.

 Then the following are equivalent:
\begin{enumerate}
\item[(1)] $K_X + \Delta$ is ample;
\item[(2)] $K_X+\Delta$ is big and its restriction to $\lfloor \Delta \rfloor$ is ample.
\end{enumerate}
If $(X, \cdel)$ is dlt, then the above conditions are also equivalent to:

\begin{enumerate}
\item[(3)] $K_X + \Delta$ is log big.

\end{enumerate}

\end{theorem}

\begin{remark}
As $(X, \Delta)$ being Mori hyperbolic implies that $K_X+\Delta$ is nef, condition $2)$ in the corollary
is equivalent to the condition stated in Theorem \ref{ample.thm}:
\[
(K_X + \Delta)^{\dim X} > 0 \; {\rm and}\; 
(K_X+\cdel)^{\dim W} \cdot W > 0, \; {\rm for \; any \;lc \; center} \; 
W.
\]
\end{remark}

\begin{remark}
The assumption on the bigness of $K_X+\Delta$ in the theorem is necessary.

In fact, the pair $(\mathbb{P}^1, \{0\} + \{\infty\})$ is log canonical and its lc centers are the points $0$ 
and $\infty$. The divisor $K_{{P}^1}+ \{0\} + \{\infty\}$ is clearly ample along the two lc centers, yet the divisor
is linearly equivalent to $0$.
\end{remark}

\begin{proof}[Proof of Theorem \ref{nakai.crit.m.hyp.cor}.]
Again, (1) implies (2) and (3).
Moreover, as $(X, \cdel)$ is Mori hyperbolic, it is nef by Theorem \ref{main.thm}.
 
Let us prove that (2) implies (1). As $K_X+\cdel$ is big and ample along $\lfloor \cdel \rfloor$, to prove its ampleness on $X$, it suffices to prove that $K_X+\cdel$ intersects all rational curves on $X$ with
strictly positive degree. Let us assume there exists a rational curve $C$ such that $(K_X+\cdel) \cdot C=0$. We can 
assume that $K_X+(\cdel - \epsilon \lfloor \cdel \rfloor) \cdot C < 0$, for any $\epsilon > 0$. Let us notice that 
$K_X+(\cdel - \epsilon \lfloor \cdel \rfloor)$ is ample along $\lfloor \cdel \rfloor$ for $0<\epsilon\ll1$. 
Passing to a dlt modification as in Theorem \ref{dlt.mod.thm}, we can assume that $X$ is $\ratl$-factorial and the proof is the same as that of Proposition \ref{mh.nonlc.dlt.prop}.

If $(3)$ holds, then by induction on $\dim X$ it follows immediately that 
$K_X+\cdel$ is ample along $\lfloor \cdel \rfloor$. Moreover, the definition of log bigness implies that $K_X+\Delta$ is also big, which terminates the proof.
\end{proof}

\bibliography{bib}

\newcommand{\etalchar}[1]{$^{#1}$}
\providecommand{\bysame}{\leavevmode\hbox to3em{\hrulefill}\thinspace}
\providecommand{\MR}{\relax\ifhmode\unskip\space\fi MR }
\providecommand{\MRhref}[2]{%
  \href{http://www.ams.org/mathscinet-getitem?mr=#1}{#2}
}
\providecommand{\href}[2]{#2}
\begin{thebibliography}{Amb05}

\bibitem[Amb03]{amb03}
F.~Ambro, \emph{Quasi-log varieties}, Tr. Mat. Inst. Steklova \textbf{240}
  (2003), no.~Biratsion. Geom. Linein. Sist. Konechno Porozhdennye Algebry,
  220--239. \MR{1993751 (2004f:14027)}

\bibitem[Amb05]{amb05}
\bysame, \emph{The moduli {$b$}-divisor of an lc-trivial fibration}, Compos.
  Math. \textbf{141} (2005), no.~2, 385--403. \MR{2134273 (2006d:14015)}

\bibitem[DCS21]{1608.02997}
Gabriele Di~Cerbo and Roberto Svaldi, \emph{Birational boundedness of
  low-dimensional elliptic {C}alabi-{Y}au varieties with a section}, Compos.
  Math. \textbf{157} (2021), no.~8, 1766--1806. \MR{4292177}

\bibitem[Dem12]{dem2012}
J.-P. Demailly, \emph{Kobayashi pseudo-metrics, entire curves and hyperbolicity
  of algebraic varieties}, 2012, Electronic document available at the webpage
  http://www-fourier.ujf-grenoble.fr/$\sim$demailly/manuscripts/kobayashi.pdf.

\bibitem[FG12]{fuj-gong12}
O.~Fujino and Y.~Gongyo, \emph{On canonical bundle formulas and
  subadjunctions}, Michigan Math. J. \textbf{61} (2012), no.~2, 255--264.
  \MR{2944479}

\bibitem[FG14]{fuj-gong.1210}
\bysame, \emph{On the moduli b-divisors of lc-trivial fibrations}, Ann. Inst.
  Fourier (Grenoble) \textbf{64} (2014), no.~4, 1721--1735. \MR{3329677}

\bibitem[Fil18]{1807.04847}
S.~Filipazzi, \emph{On a generalized canonical bundle formula and generalized
  adjunction}, 2018, ArXiv e-print, arXiv:1807.04847v3, to appear in Annali
  della Scuola Normale Superiore di Pisa - Classe di Scienze.

\bibitem[FS20]{mio.connectedness}
S.~Filipazzi and R.~Svaldi, \emph{On the connectedness principle and dual
  complexes for generalized pairs}, 2020, Arxiv e-print, arXiv:2010.08018v2.

\bibitem[Fuj09]{fuj09}
O.~Fujino, \emph{Introduction to the log minimal model program for log
  canonical pairs}, ArXiv e-print (2009), Arxiv:0907.1506v1.

\bibitem[Fuj11]{fuj.nonvan.12}
\bysame, \emph{Non-vanishing theorem for log canonical pairs}, J. Algebraic
  Geom. \textbf{20} (2011), no.~4, 771--783. \MR{2819675 (2012j:14021)}

\bibitem[HM07]{hac-mck05}
C.~D. Hacon and J.~M\textsuperscript{c}Kernan, \emph{On {S}hokurov's rational
  connectedness conjecture}, Duke Math. J. \textbf{138} (2007), no.~1,
  119--136. \MR{2309156 (2008f:14030)}

\bibitem[K{\etalchar{+}}92]{koll-etal-book}
J.~Koll{\'a}r et~al., \emph{Flips and abundance for algebraic threefolds},
  Flips and abundance for algebraic threefolds, Ast{\'e}risque, vol. 211, SMF,
  1992, pp.~i+258.

\bibitem[Kaw91]{kawam91}
Y.~Kawamata, \emph{On the length of an extremal rational curve}, Invent. Math.
  \textbf{105} (1991), no.~3, 609--611. \MR{1117153 (92m:14026)}

\bibitem[KK10]{koll-kov10}
J.~Koll{\'a}r and S.~J. Kov{\'a}cs, \emph{Log canonical singularities are {D}u
  {B}ois}, J. Amer. Math. Soc. \textbf{23} (2010), no.~3, 791--813. \MR{2629988
  (2011m:14061)}

\bibitem[KM98]{koll-mor-book}
J.~Koll{\'a}r and S.~Mori, \emph{Birational geometry of algebraic varieties},
  Cambridge Tracts in Mathematics, vol. 134, Cambridge University Press,
  Cambridge, 1998, With the collaboration of C. H. Clemens and A. Corti,
  Translated from the 1998 Japanese original. \MR{1658959 (2000b:14018)}

\bibitem[Kol13]{koll-book13}
J\'{a}nos Koll\'{a}r, \emph{Singularities of the minimal model program},
  Cambridge Tracts in Mathematics, vol. 200, Cambridge University Press,
  Cambridge, 2013, With a collaboration of S\'{a}ndor Kov\'{a}cs. \MR{3057950}

\bibitem[Laz04a]{laz1.book}
R.~Lazarsfeld, \emph{Positivity in algebraic geometry. {I}}, Ergebnisse der
  Mathematik und ihrer Grenzgebiete. 3. Folge. A Series of Modern Surveys in
  Mathematics [Results in Mathematics and Related Areas. 3rd Series. A Series
  of Modern Surveys in Mathematics], vol.~48, Springer-Verlag, Berlin, 2004,
  Classical setting: line bundles and linear series. \MR{2095471
  (2005k:14001a)}

\bibitem[Laz04b]{laz2.book}
\bysame, \emph{Positivity in algebraic geometry. {II}}, Ergebnisse der
  Mathematik und ihrer Grenzgebiete. 3. Folge. A Series of Modern Surveys in
  Mathematics [Results in Mathematics and Related Areas. 3rd Series. A Series
  of Modern Surveys in Mathematics], vol.~49, Springer-Verlag, Berlin, 2004,
  Positivity for vector bundles, and multiplier ideals. \MR{2095472
  (2005k:14001b)}

\bibitem[LZ17]{lu-zh}
Steven S.~Y. Lu and De-Qi Zhang, \emph{Positivity criteria for log canonical
  divisors and hyperbolicity}, J. Reine Angew. Math. \textbf{726} (2017),
  173--186. \MR{3641656}

\bibitem[Mat02]{matsuki.book}
K.~Matsuki, \emph{Introduction to the {M}ori program}, Universitext,
  Springer-Verlag, New York, 2002. \MR{1875410 (2002m:14011)}

\bibitem[MP12]{mcq-pac}
M.~McQuillan and G.~Pacienza, \emph{Remarks about bubbles}, Arxiv e-print
  (2012), Arxiv:1211.0203v1.

\bibitem[OX12]{MR2955764}
Y.~Odaka and C.~Xu, \emph{Log-canonical models of singular pairs and its
  applications}, Math. Res. Lett. \textbf{19} (2012), no.~2, 325--334.
  \MR{2955764}

\bibitem[SS19]{mio.term}
C.~Spicer and R.~Svaldi, \emph{Local and global applications of the minimal
  model program for co-rank one foliations on threefolds}, 2019, ArXiv e-print,
  arXiv:1908.05037v1.

\end{thebibliography}

\end{document}